\documentclass[12pt,a4paper,twoside]{amsart}

\usepackage{mathrsfs}          
\usepackage[all]{xy}    
\usepackage[latin1]{inputenc}
\usepackage{amscd}
\usepackage{latexsym}
\usepackage{multicol}
\usepackage{amsmath}
\usepackage{amssymb}
\usepackage{amsthm}

\theoremstyle{plain}
\newtheorem{Theorem}{Theorem}[section]

\newtheorem{corollary}[Theorem]{Corollary}
\newtheorem{proposition}[Theorem]{Proposition}

\newtheorem{conjecture2*}{Conjecture}
\newtheorem{remark2*}[conjecture2*]{Remark}
\newtheorem{Theorem2*}[conjecture2*]{Theorem}
\newtheorem{proposition2*}[conjecture2*]{Proposition}
\newtheorem{claim2*}[conjecture2*]{Claim}
\newtheorem{example2*}[conjecture2*]{Example}
\newtheorem{lemma2*}[conjecture2*]{Lemma}

\theoremstyle{definition}

\newtheorem{example}[Theorem]{Example}

\newtheorem{remark}[Theorem]{Remark}

\newtheorem{pkt}[Theorem]{}
\newtheorem{construction}[Theorem]{Construction}

\newcommand{\sL}{{\mathcal L}}
\newcommand{\sM}{{\mathcal M}}

\newcommand{\sO}{{\mathcal O}}

\newcommand{\sY}{{\mathcal Y}}

\newcommand{\A}{{\mathbb A}}

\newcommand{\C}{{\mathbb C}}

\newcommand{\E}{{\mathbb E}}
\newcommand{\F}{{\mathbb F}}

\newcommand{\N}{{\mathbb N}}
\newcommand{\bP}{{\mathbb P}}

\title{Examples of Calabi-Yau 3-manifolds with complex multiplication}
\author{Jan Christian Rohde}\thanks{This paper has been financially supported by the Leibniz-Preis awarded to H\'el\`ene Esnault and
Eckart Viehweg.}
\email{jan.rohde@uni-due.de}
\begin{document}
\maketitle

\section*{Introduction}

By string theoretical considerations, one is interested in Calabi-Yau manifolds, since Calabi-Yau
3-manifolds provide conformal field theories (CFT). One is especially interested in Calabi-Yau
3-manifolds with complex multiplication, since such a manifold has many symmetries and mirror
pairs of Calabi-Yau 3-manifolds with complex multiplication yield rational conformal field
theories (RCFT) (see \cite{GV}).

E. Viehweg and K. Zuo \cite{VZ5} have constructed a family of Calabi-Yau
$3$-manifolds with dense set of complex multiplication fibers. This construction is given by a
tower of cyclic coverings. We will use a similar construction to obtain $K3$ surfaces with complex
multiplication. 

C. Voisin \cite{Voi2} has
described a method to obtain Calabi-Yau 3-manifolds by using
involutions on $K3$ surfaces. C. Borcea \cite{Bc2} has independently
arrived at a more general version of the latter method, which allows to
construct Calabi-Yau manifolds in arbitrary
dimension. By using this method and $K3$ surfaces with complex
multiplication, we will obtain our concrete examples of Calabi-Yau 3-manifolds with complex
multiplication.

Our methods are very similar to the methods in \cite{JCR}, which contains different
concrete examples of Calabi-Yau 3-manifolds with complex multiplication. The examples of \cite{JCR}
occur as fibers of a family with a dense set of complex multiplication fibers. Here we give some
examples, which are not necessarily fibers of a non-trivial family with a dense set of complex
multiplication fibers. The first two sections give two different classes of examples by using
involutions on $K3$ surfaces.

In the third section we will prove that a $K3$ surface with a degree 3 automorphism has complex
multiplication. By using methods, which has been introduced in \cite{JCR} Section $9.1$ and
Section $9.2$, we will use this automorphism and the Fermat curve of degree 3 for the construction
of a Calabi-Yau 3-manifold with complex multiplication.

We use the same methods as in \cite{JCR} Chapter 10 to determine the Hodge numbers of our
examples.

\tableofcontents

\section{Construction by degree 2 coverings of a ruled surface}
We use similar methods as in \cite{JCR}. Hence we start by finding curves with complex
multiplication. The following theorem yields some examples:

\begin{Theorem} \label{geilomat}
Let $0 <d_1, d < m$, and $\xi_k$ denote a primitive $k$-th. root of
unity for all $k \in \N$. Then  the curve $C$, which is locally given by
$$y^m = x^{d_1}\prod\limits_{i = 1}^{n-2} (x-\xi_{n-2}^i)^d,$$
is covered by the Fermat curve $\F_{(n-2)m}$ locally given by
$$y^{(n-2)m}+x^{(n-2)m}+1=0$$
and has complex multiplication.
\end{Theorem}
\begin{proof}
(see \cite{JCR} Theorem $2.4.4$)
\end{proof}

\begin{example} \label{exa1}
By the preceding theorem, the curves locally given by
$$y^4 = x_1^8 + x_0^8, \ \ y^4 = x_1(x_1^7 + x_0^7), \ \ y^4 = x_1(x_1^6 + x_0^6)x_0$$
have complex multiplication. These curves are degree 4 covers of the projective line and have the
genus 9 as one can easily calculate by the Hurwitz formula.
\end{example}

The curves of the preceding example have a natural interpretation as cyclic covers of $\bP^1$ of
degree 4. One can identify these covers with the set of their 8 branch points in $\bP^1$. Thus let
$\sM_8$ denote the configuration space of 8 different points in $\bP^1$.\footnote{Note that this
is not the same notation as in \cite{JCR}.} We use a modified version of the construction in
\cite{VZ5}, Section 5 to construct $K3$ surfaces with complex multiplication by Example
$\ref{exa1}$ in a first step. This method is nearly the same method as in \cite{JCR} Section $8.2$.

For our application, it is sufficient to work with $\bP^1$-bundles
over $\bP^1$ resp., with rational ruled surfaces. Let $\pi_n: \bP_n \to \bP^1$ denote
the rational ruled surface given by $\bP(\sO_{\bP^1} \oplus
\sO_{\bP^1}(n))$ and $\sigma$ denote a non-trivial global section
of $\sO_{\bP^1}(8)$, which has the 8 different zero points represented by a point $q \in \sM_8$.
The sections $E_{\sigma}$, $E_0$ and $E_{\infty}$ of $\bP(\sO \oplus \sO(8))$ are
induced by
$${\rm id}\oplus \sigma:\sO \to \sO\oplus \sO(8), \ \ 
{\rm id}\oplus 0:\sO \to \sO\oplus \sO(8)$$
$$\mbox{and} \  \ 0\oplus {\rm id}:\sO(8) \to \sO\oplus \sO(8)$$
resp., by the corresponding surjections onto the cokernels of these embeddings as described in
\cite{hart}, {\bf II}. Proposition 7.12.

\begin{remark}
The divisors $E_{\sigma}$ and $E_0$ intersect each other transversally over the 8 zero points of
$\sigma$. Recall that ${\rm Pic}(\bP_8)$ has a basis given by a fiber and an arbitrary section.
Hence by the fact that $E_{\sigma}$ and $E_0$ do not intersect $E_{\infty}$, one concludes that
they are linearly equivalent with self-intersection number 8. Since $E_{\infty}$
is a section, it intersects each fiber transversally. Thus one has that
$E_{\infty} \sim  E_0-(E_0.E_0)F$, where $F$ denotes a fiber. Hencefore one concludes
$$E_{\infty}.E_{\infty} = E_{\infty}.(E_0-(E_0.E_0)F) = -(E_0.E_0) = -8.$$
\end{remark}

Next we establish a morphism $\mu:\bP_2 \to \bP_8$ over $\bP^1$. By
\cite{hart}, ${\bf II}$. Proposition 7.12., this is the same as to
give a surjection $\pi_2^*(\sO \oplus \sO(8)) \to \sL$, where $\sL$
is an invertible sheaf on $\bP_2$. By the composition
$$\pi_2^*(\sO \oplus \sO(8)) = \pi_2^*(\sO) \oplus \pi_2^*\sO(8) \hookrightarrow
\bigoplus\limits_{i = 0}^4 \pi_2^*\sO(2i) = Sym^4(\pi_2^*(\sO \oplus \sO(8))) \to
\sO_{\bP_2}(4),$$
where the last morphism is induced by the natural surjection $\pi_2^*(\sO \oplus \sO(2)) \to
\sO_{\bP_2}(1)$ (see \cite{hart}, {\bf II}. Proposition 7.11), we obtain a morphism $\mu^*$ of
sheaves. This morphism $\mu^*$ is not a surjection onto $\sO_{\bP_2}(4)$, but onto its image
$\sL \subset \sO_{\bP_2}(4)$. Over $\A^1 \subset \bP^1$ all rational ruled surfaces are locally
given by ${\rm Proj}(\C[x])[y_1,y_2]$, where $x$ has the weight 0. Hence we have locally that
$\pi_2^*(\sO \oplus \sO(8)) = \sO e_1\oplus \sO e_2$. Over $\A^1$ the morphism $\mu^*$ is given by
$$e_1 \to y_1^{4}, e_2 \to y_2^{4}$$
such that the sheaf $\sL = im(\mu^*) \subset \sO_{\bP_2}(4)$ is invertible.
Thus the morphism $\mu: \bP_2 \to \bP_8$ corresponding to $\mu^*$ is locally
given by the ring homomorphism
$$(\C[x])[y_1,y_2] \to (\C[x])[y_1,y_2] \ \ \mbox{via} \  \ y_1 \to y_1^8 \ \ \mbox{and}
\ \ y_2 \to y_2^8.$$

\begin{construction}\label{4.1}
One has a commutative diagram
$$\xymatrix{
{\sY'} \ar[rr]^{\tau'} &   & {\bP'_2} \ar[rr]^{\mu'} &   & {\bP^1 \times \bP^1}\\
{\hat \sY} \ar[rr]^{\hat \tau} \ar[d]_{\rho} \ar[u]^{\delta} &  & {\hat \bP_2} \ar[rr]^{\hat \mu}
\ar[d]^{\rho_2} \ar[u]_{\delta_2} &  & {\hat \bP_8} \ar[d]^{\rho_8} \ar[u]_{\delta_8}\\
{\sY} \ar[rr]^{\tau}_{\sqrt[2]{\frac{\mu^*E_\sigma}{3\cdot (\mu^*E_0)_{red}}}}
\ar[d]_{\pi} &  & {\bP_2} \ar[rr]^{\mu}_{\sqrt[4]{\frac{E_\infty + 8 \cdot F}{E_0}}} \ar[d]^{\pi_2}
&  & {\bP_8} \ar[d]^{\pi_8}\\
{\bP^1} \ar[rr]^{\rm id} &   & {\bP^1} \ar[rr]^{\rm id} &   & {\bP^1}
}$$
of morphisms between normal varieties with:
\begin{enumerate}
\item[(a)] $\delta$, $\delta_2$, $\delta_8$, $\rho$, $\rho_2$ and $\rho_8$
are birational.
\item[(b)] $\pi$ is a family of curves, $\pi_2$ and $\pi_8$
are $\bP^1$-bundles.
\end{enumerate}
\end{construction}
\begin{proof}
One must only explain $\delta_8$ and $\rho_8$. Recall that
$E_{\sigma}$ is a section of $\bP(\sO \oplus \sO(8))$, which
intersects $E_0$ transversally in exactly 8 points. The morphism $\rho_8$ is the
blowing up of the 8 intersection points of $E_0 \cap E_{\sigma}$.
The preimage of the 8 points given by $q \in \sM_8$ with respect to $\pi_8
\circ \rho_8$ consists of the exceptional divisor $\hat D_1$ and the
proper transform $\tilde D_2$ of the preimage of these 8 points
with respect to $\rho_8$ given by 8 rational curves with self-intersection number $-1$. The
morphism $\delta_8$ is obtained by blowing down $\tilde D_2$.
\end{proof}

\begin{remark}
The section $\sigma$ has the zero divisor given by some $q \in \sM_8$. Hence one obtains
$\mu^*(E_{\sigma}) \cong C$, where $C \to \bP^1$ is the cyclic cover of degree 4 as in Example $\ref{exa1}$
ramified over the 8 points given by $\sigma$. The surface $\sY$ is a cyclic degree 2 cover of $\bP_2$
ramified over $C$. Thus it has an involution. It is given by the invertible sheaf
$$\sL = \omega_{\bP_2}^{-1}$$
and the divisor
$$B = \mu^*(E_{\sigma}), \ \ \mbox{where} \ \ \sO(B) \cong \sL^{2},$$
with the notation of \cite{Bart} {\bf I}. 17. Thus \cite{Bart} {\bf I}. Lemma $17.1$ implies that
$\sY$ is a $K3$ surface.
\end{remark}

By \cite{JCR} Lemma $10.4.1$, there is only one elliptic curve with a cyclic degree 4 cover onto
$\bP^1$. Let $\E$ denote this curve, which is locally given by
$$y^4 = x(x-1)^2.$$
One can easily see that $\E$ has the $j$ invariant 1728. Thus
$\E$ has complex multiplication.

We introduce a new notation. Let $n \in \N$, let $\xi$ be a fixed primitive $n$-th. root of unity
and let $C_1$ and $C_2$ be curves locally given by
$$y^n=f_1(x) \  \ \mbox{and} \  \ y^n=f_2(x),$$
where $f_1, f_2 \in \C[x]$. By $(x,y) \to (x,\xi y)$, one can define an automorphism $\gamma_i$
on $C_i$ for $i =1,2$. The surface $C_1 \times C_2/\langle(1, 1)\rangle$ is the quotient of
$C_1 \times C_2$ by $\langle (\gamma_1,\gamma_2)\rangle$.

\begin{proposition} \label{proposition}
The surface $\sY$ is birationally equivalent to
$C \times \E/\langle(1, 1)\rangle$.\footnote{Similarly to
\cite{VZ5}, Construction 5.2, we show that $\sY'$ is birationally equivalent to
$C \times \E/\langle(1, 1)\rangle$.}
\end{proposition}
\begin{proof}
Let $\tilde  E_{\bullet}$ denote the proper transform of the section $E_{\bullet }$ with respect
to $\rho_8$. Then $\hat \mu$ is the Kummer covering given by
$$\sqrt[4]{\frac{\tilde E_\infty + 8 \cdot F}{\tilde E_0+ \hat D_1}},$$
where $\hat D_1$ denotes the exceptional divisor of $\rho_8$. Thus the morphism $\mu'$ is the
Kummer covering
$$\sqrt[4]{\frac{ (\delta_8)_* \tilde E_\infty + 8 \cdot (\delta_8)_* F}
{(\delta_8)_*\tilde E_0+ (\delta_8)_*\hat D_1}}
=\sqrt[4]{\frac{ \bP^1 \times \{\infty\} + 8 \cdot (P \times \bP^1)}
{\bP^1 \times \{0\} + \Delta  \times\bP^1}},$$
where $\Delta $ is the divisor of the 8 different points in $\bP^1$ given by $q \in \sM_8$
and $P \in \bP^1$ is the point with the fiber $F$. Since $E_0 + E_{\sigma}$ is a normal crossing
divisor, $\tilde E_{\sigma}$ neither meets $\tilde E_{0}$
nor $\tilde D_2$, where $\tilde D_2$ is the proper transform of $\pi_8^*(\Delta)$. Therefore
$(\delta_8)_*\tilde E_{\sigma }$ neither meets
$$(\delta_8)_*\tilde E_0 = \bP^1 \times \{0\} \  \ \mbox{nor} \  \
(\delta_8)_*\tilde E_{\infty} = \bP^1 \times \{\infty\}.$$
Hence one can choose coordinates in $\bP^1$ such that
$(\delta_8)_*\tilde E_{\sigma} = \bP^1 \times \{1\}$.

By the definition of $\tau$, we obtain that $\hat \tau$ is given by
$$\sqrt[2]{\frac{ \rho_2^*\mu^*(E_{\sigma} )}{\rho_2^*\mu^*(E_0)}}
= \sqrt[2]{\frac{ \hat \mu^*(\tilde E_{\sigma} )}{\hat \mu^*(\tilde E_0)}},$$
and $\tau'$ is given by
$$\sqrt[2]{\frac{ \mu'^*(\bP^1 \times \{1\})}
{\mu'^*(\bP^1 \times \{0\})}}.$$
By the fact that the last function is the root of the pullback of a function
on $\bP^1 \times \bP^1$ with respect to $\mu'$, it is possible to reverse the order of the field
extensions corresponding to $\tau'$ and $\mu'$ such that the resulting varieties obtained by
Kummer coverings are birationally equivalent. Hence we have the composition of $\beta:
\bP^1 \times \bP^1 \to \bP^1 \times \bP^1$ given by
$$\sqrt[2]{\frac{ \bP^1 \times \{1\}}{\bP^1 \times \{0\}}}$$
with
$$\sqrt[4]{\frac{ \beta^*(\bP^1 \times \{\infty\}) + 8 \cdot  (P \times \bP^1)}
{\beta^*(\bP^1 \times \{0\}) + (\Delta \times\bP^1)}},$$
which yields the covering variety isomorphic to $\E \times C/\langle(1,1)\rangle$.
\end{proof}

As in \cite{JCR} Section $8.2$ we conclude:

\begin{corollary} \label{cmcm}
If the curve $C$ has complex multiplication, the $K3$-surface $\sY$ has only
commutative Hodge groups.
\end{corollary}

By the the preceding corollary, our Example $\ref{exa1}$ yields 3 $K3$ surfaces with complex multiplication
locally given by
$$y_2^2+y_1^4 + x_1^8 + x_0^8, \ \ y_2^2+y_1^4 + x_1(x_1^7 + x_0^7), \ \
y_2^2 + y_1^4 + x_1(x_1^6 + x_0^6)x_0.$$

\begin{proposition} \label{prop}
For $i = 1,2$ assume that $C_i$ is a Calabi-Yau $i$-manifold with complex multiplication
endowed with the involution $\iota_i$ such that $\iota_i$ acts by $-1$ on $\Gamma(\omega_{C_i})$.
By blowing up the singular locus of $C_1 \times C_2/\langle(\iota_1, \iota_2)\rangle$,
one obtains a Calabi-Yau $3$-manifold with complex multiplication.
\end{proposition}
\begin{proof}
It is well-known that an involution on a Calabi-Yau 2-manifold resp., a $K3$ surface, which acts by
$-1$ on $\Gamma(\omega)$, has a smooth divisor of fixed points or it has not any fixed point. Thus
the proof follows from \cite{JCR} Section $7.2$.
\end{proof}

Now we need some elliptic curves with complex multiplication:

\begin{example} \label{exa2}
Elliptic curves with $CM$ has been well studied by number theorists. Some examples of elliptic curves
with complex multiplication are given by the following list:
\begin{center}
\begin{tabular}{|c|c|} \hline
equation & $j$ invariant \\ \hline \hline
$y_1^2x_0 = x_1^3-x_0^3$ & 0\\ \hline
$y_1^2x_0 = x_1(x_1-x_0)(x_1-2x_0)$ & 1728 \\ \hline
$y^2x_0 = x_1(x_1-x_0)(x_1-(1+\sqrt{2})^2x_0)$ & 8000 \\ \hline
$y^2x_0 = x_1(x_1-x_0)(x_1-\frac{1}{4}(3+i\sqrt{7})^2x_0)$ & $-3375$ \\ \hline
$y^2x_0 = x_1^3 - 15x_1x_0^2 +22x_0^3$ & 54000 \\ \hline
$y^2x_0 = x^3 -595x_1x_0^2+5586x_0^3$ & 16581375 \\ \hline
\end{tabular}
\end{center}
Note that the equations allow an explicit definition of an involution on these elliptic curves.
(see \cite{JCR} Section $7.4$)
\end{example}

\begin{pkt} \label{punkt}
By combining our 3 examples of $K3$ surfaces and the 6 elliptic curves and using Propostion $\ref{prop}$, we
have 18 examples of Calabi-Yau 3-manifolds with complex multiplication. It seems to be quite easy
to describe these examples by explicit equations. By \cite{Voi2}, one has equations to determine
the Hodge numbers of these examples. Let $C_2$ be a $K3$ surface satisfying the assumptions of
Proposition $\ref{prop}$, let $N$ be the number of curves in the ramification locus of the quotient
map $C_2 \to C_2/\iota_2$ and let $N'$ be given by
$$N' = g_1 + \ldots + g_N,$$
where $g_i$ denotes the genus of the $i$-th. curve in the ramification locus. Then one has for the
Calabi-Yau 3-manifold, which results by Proposition $\ref{prop}$:
$$h^{1,1} = 11+5N-N'$$
$$h^{2,1} = 11+5N'-N$$
Thus in our case the Hodge numbers are given by
$$h^{1,1} =7 \  \ \mbox{and} \  \ h^{2,1} = 55.$$
\end{pkt}

\section{Construction by degree 2 coverings of $\bP^2$}

\begin{example} \label{exa3}
By Theorem $\ref{geilomat}$, the projective curves given by
$$y^6 = x_1^6 + x_0^6, \ \ y^6 = x_1(x_1^5 + x_0^5), \ \ y^6 = x_1(x_1^4 + x_0^4)x_0$$
have complex multiplication. These curves have the genus 10 as one can easily calculate by the
Hurwitz formula.
\end{example}

Let $\sM_6$ denote the configuration space of 6 different points in $\bP^1$. Again we use a modified
version of the construction in \cite{VZ5}, Section 5.

Here the sections $E_{\sigma}$, $E_0$ and $E_{\infty}$ of $\bP(\sO \oplus \sO(6))$ are induced by
$${\rm id}\oplus \sigma:\sO \to \sO\oplus \sO(6), \ \ 
{\rm id}\oplus 0:\sO \to \sO\oplus \sO(6)$$
$$\mbox{and} \  \ 0\oplus {\rm id}:\sO(6) \to \sO\oplus \sO(6)$$
resp., by the corresponding surjections onto the cokernels of these embeddings as described in
\cite{hart}, {\bf II}. Proposition 7.12.

One concludes similarly to the preceding section that
$$E_{\infty}.E_{\infty} = E_{\infty}.(E_0-(E_0.E_0)F) = -(E_0.E_0) = -6.$$

By the composition
$$\pi_1^*(\sO \oplus \sO(6)) = \pi_1^*(\sO) \oplus \pi_1^*\sO(6) \hookrightarrow
\bigoplus\limits_{i = 0}^6 \pi_1^*\sO(i) = Sym^6(\pi_1^*(\sO \oplus \sO(6))) \to
\sO_{\bP_1}(6),$$
where the last morphism is induced by the natural surjection $\pi_2^*(\sO \oplus \sO(1)) \to
\sO_{\bP_1}(1)$ (see \cite{hart}, {\bf II}. Proposition 7.11), we obtain a morphism $\mu^*$ of
sheaves as in the preceding section. 
The morphism $\mu: \bP_1 \to \bP_1$ corresponding to $\mu^*$ is locally
given by the ring homomorphism
$$(\C[x])[y_1,y_2] \to (\C[x])[y_1,y_2] \ \ \mbox{via} \  \ y_1 \to y_1^6 \ \ \mbox{and}
\ \ y_2 \to y_2^6.$$

\begin{construction}
One has a commutative diagram
$$\xymatrix{
{\sY'} \ar[rr]^{\tau'} &   & {\bP'_1} \ar[rr]^{\mu'} &   & {\bP^1 \times \bP^1}\\
{\hat \sY} \ar[rr]^{\hat \tau} \ar[d]_{\rho} \ar[u]^{\delta} &  & {\hat \bP_1} \ar[rr]^{\hat \mu}
\ar[d]^{\rho_1} \ar[u]_{\delta_1} &  & {\hat \bP_6} \ar[d]^{\rho_6} \ar[u]_{\delta_6}\\
{\sY} \ar[rr]^{\tau}_{\sqrt[2]{\frac{\mu^*E_\sigma}{3\cdot (\mu^*E_0)_{red}}}}
\ar[d]_{\pi} &  & {\bP_1} \ar[rr]^{\mu}_{\sqrt[6]{\frac{E_\infty + 6 \cdot F}{E_0}}} \ar[d]^{\pi_1}
&  & {\bP_6} \ar[d]^{\pi_6}\\
{\bP^1} \ar[rr]^{\rm id} &   & {\bP^1} \ar[rr]^{\rm id} &   & {\bP^1}
}$$
of morphisms between normal varieties with:
\begin{enumerate}
\item[(a)] $\delta$, $\delta_1$, $\delta_1$, $\rho$, $\rho_1$ and $\rho_6$
are birational.
\item[(b)] $\pi$ is a family of curves, $\pi_1$ and $\pi_6$
are $\bP^1$-bundles.
\end{enumerate}
\end{construction}
\begin{proof}
One must only explain $\delta_6$ and $\rho_6$. These morphisms are given by blowing up morphisms
similar to the preceding section.
\end{proof}

\begin{remark}
The section $\sigma$ has the zero divisor given by some $q \in \sM_6$. Hence one obtains
$\mu^*(E_{\sigma}) \cong C$, where $C \to \bP^1$ is a cyclic cover of degree 6 as in Example $\ref{exa3}$
ramified over the 6 points given by $\sigma$. The surface $\sY$ is a cyclic degree 2 cover of $\bP_1$
ramified over $C$. Thus it is birationally equivalent to the $K3$ surface given the degree 2 cover
of $\bP^2$ ramified over $C$.
\end{remark}

Let $C'$ denote the projective smooth curve locally given by
$$y^6 = x(x-1).$$
By Theorem $\ref{geilomat}$, it has complex multiplication.

\begin{proposition}
The surface $\sY$ is birationally equivalent to
$C \times C'/\langle(1, 1)\rangle$.
\end{proposition}
\begin{proof}
Let $\tilde  E_{\bullet}$ denote the proper transform of the section $E_{\bullet }$ with respect
to $\rho_6$. Then $\hat \mu$ is the Kummer covering given by
$$\sqrt[6]{\frac{\tilde E_\infty + 6 \cdot F}{\tilde E_0+ \hat D_1}},$$
where $\hat D_1$ denotes the exceptional divisor of $\rho_6$. Thus the morphism $\mu'$ is the
Kummer covering
$$\sqrt[6]{\frac{ (\delta_6)_* \tilde E_\infty + 6 \cdot (\delta_6)_* F}
{(\delta_6)_*\tilde E_0+ (\delta_6)_*\hat D_1}}
=\sqrt[6]{\frac{ \bP^1 \times \{\infty\} + 6 \cdot (P \times \bP^1)}
{\bP^1 \times \{0\} + \Delta  \times\bP^1}},$$
where $\Delta $ is the divisor of the 6 different points in $\bP^1$ given by $q \in \sM_6$
and $P \in \bP^1$ is the point with the fiber $F$. Since $E_0 + E_{\sigma}$ is a normal crossing
divisor, $\tilde E_{\sigma}$ neither meets $\tilde E_{0}$
nor $\tilde D_2$, where $\tilde D_2$ is the proper transform of $\pi_6^*(\Delta)$. Therefore
$(\delta_6)_*\tilde E_{\sigma }$ neither meets
$$(\delta_6)_*\tilde E_0 = \bP^1 \times \{0\} \  \ \mbox{nor} \  \
(\delta_6)_*\tilde E_{\infty} = \bP^1 \times \{\infty\}.$$
Hence one can choose coordinates in $\bP^1$ such that
$(\delta_6)_*\tilde E_{\sigma} = \bP^1 \times \{1\}$.

By the definition of $\tau$, we obtain that $\hat \tau$ is given by
$$\sqrt[2]{\frac{ \rho_1^*\mu^*(E_{\sigma} )}{\rho_1^*\mu^*(E_0)}}
= \sqrt[2]{\frac{ \hat \mu^*(\tilde E_{\sigma} )}{\hat \mu^*(\tilde E_0)}},$$
and $\tau'$ is given by
$$\sqrt[2]{\frac{ \mu'^*(\bP^1 \times \{1\})}
{\mu'^*(\bP^1 \times \{0\})}}.$$
By the fact that the last function is the root of the pullback of a function
on $\bP^1 \times \bP^1$ with respect to $\mu'$, it is possible to reverse the order of the field
extensions corresponding to $\tau'$ and $\mu'$ such that the resulting varieties obtained by
Kummer coverings are birationally equivalent. Hence we have the composition of $\beta:
\bP^1 \times \bP^1 \to \bP^1 \times \bP^1$ given by
$$\sqrt[2]{\frac{ \bP^1 \times \{1\}}{\bP^1 \times \{0\}}}$$
with
$$\sqrt[6]{\frac{ \beta^*(\bP^1 \times \{\infty\}) + 6\cdot  (P \times \bP^1)}
{\beta^*(\bP^1 \times \{0\}) + (\Delta \times\bP^1)}},$$
which yields the covering variety isomorphic to $C' \times C/\langle(1,1)\rangle$.
\end{proof}

Hence $\sY$ is birationally equivalent to $C' \times C/\langle(1,1)\rangle$. As in \cite{JCR}
Section $8.2$ we conclude:

\begin{corollary}
If the curve $C$ has complex multiplication, the $K3$-surface $\sY$ has only
commutative Hodge groups.
\end{corollary}

\begin{pkt}
By the preceding corollary, our Example $\ref{exa3}$ yields 3 different $K3$ surfaces with complex multiplication
as degree 2 covers of $\bP^2$, which are locally given by
$$y_2^2+ y_1^6 = x_1^6 + x_0^6, \ \ y_2^2+ y_1^6 = x_1(x_1^5 + x_0^5),
\ \ y_2^2+y_1^6 = x_1(x_1^4 + x_0^4)x_0.$$
By an elliptic curve with complex multiplication, these $K3$ surfaces yield Calabi-Yau 3-manifolds
with complex multiplication. We obtain 18 Calabi-Yau 3-manifolds with complex multiplication by
using Example $\ref{exa2}$. By the same methods as in $\ref{punkt}$, one
calculates easily that the resulting Calabi-Yau 3-manifolds have the Hodge
numbers
$$h^{1,1} = 6 \ \ \mbox{and} \ \ h^{2,1} = 60$$
\end{pkt}

\section{Construction by a degree 3 quotient}

Consider the $K3$ surface
$$S = V((y_2^3-y_1^3)y_1+(x_1^3-x_0^3)x_0) \subset \bP^3.$$
By using the partial derivatives of the defining equation, one can easily verify that $S$ is smooth.
First we will prove that this surface has complex multiplication. In a second step we consider an
automorphism of degree 3 on this surface, which allows the construction of a Calabi-Yau 3-manifold
with complex multiplication.

\begin{proposition} \label{lupsi}
The $K3$ surface $S$ has complex multiplication.
\end{proposition}
\begin{proof}
Consider the isomorphic curves
$$C_1 = V(z_1^4-(y_2^3-y_1^3)y_1) \subset \bP^2,$$
$$C_2 = V(z_2^4-(x_1^3-x_0^3)x_0) \subset \bP^2.$$
Since the elliptic curve with $j$ invariant 0 given by
$$V(y^2x_0+x_1^3+x_0^3) \subset \bP^2$$ has complex multiplication,
one concludes as in \cite{JCR} Remark $7.4.2$ that $C_1$ and $C_2$ have complex multiplication, too.
The $K3$ surface $S$ is birationally equivalent to
$C_1 \times C_2/\langle(1,1)\rangle$. This follows from the rational map $C_1 \times C_2 \to S$
given by
$$((z_1:y_2:y_1),(z_2:x_1:x_0)) \to (\frac{z_2}{z_1}y_2:\frac{z_2}{z_1}y_1:x_1:x_0).$$
Thus $S$ has complex multiplication.\footnote{In \cite{Bc2} Section 5 one finds a similar rational map}
\end{proof}

\begin{pkt}
Let $\xi$ denote $e^{\frac{2\pi i}{3}}$. The $K3$ surface $S$ has an automorphism
$\gamma$ of degree 3 given by
$$(y_2:y_1:x_1:x_0) \to (\xi y_2:y_1:\xi x_1:x_0).$$
On $\{x_0 = 1\}$ we have the 4 fixed points given by 
$$(0:\sqrt[4]{-1}:0:1).$$
Now assume $x_0 = 0$. This yields
$$(y_2^3-y_1^3)y_1 = 0.$$
Thus in addition the line given by $y_1 = x_0 = 0$ is fixed.
\end{pkt}

\begin{proposition}
The automorphism $\gamma$ acts via pullback by $\xi^2$ on $\Gamma(\omega_S)$.
\end{proposition}
\begin{proof}
The $-1$ eigenspace in $\Gamma(\omega_{C_1})$ and $\Gamma(\omega_{C_2})$ comes from the cohomology
of the elliptic curve $E_0$ given by
$$y^2x_0=x_1^3-x_0^3$$
(see \cite{JCR} Section $4.2$). By the rational map in the proof of Proposition $\ref{lupsi}$, one
concludes that $\Gamma(\omega_S)$ is given by tensor product of the $-1$ eigenspace in
$\Gamma(\omega_{C_1})$ and $\Gamma(\omega_{C_2})$.

The automorphism $\gamma_{\F_3}:E_0 \to E_0$ given by $x_1 \to \xi x_1$ is the generator of the Galois group
of the degree 3 cover, which allows an identification of $E_0$ with the Fermat curve $\F_3$ of
degree 3. It acts via pullback by $\xi$ on $\Gamma(\omega_{\F_3})$. Thus the corresponding
automorphisms $\varphi_{C_1}$ and $\varphi_{C_2}$ act by $\xi$ on the $-1$ eigenspace with respect to
$C_1$ and $C_2$. Note that $(\varphi_{C_1},\varphi_{C_2})$ yields an automorphism of $(C_1,C_2)/\langle
(1,1)\rangle$. By the birational map to $S$, this automorphism corresponds to $\gamma$ and one verifies
easily that $\gamma$ acts via pullback by $\xi^2$ on $\Gamma(\omega_S)$.
\end{proof}

\begin{pkt}
Consider the blowing up $\tilde \bP^3$ of $\bP^3$ with respect to $\{y_2 = x_1 = 0\}$. Let $\tilde S$ denote the
proper transform of blowing up of $S$ with respect to the latter blowing up, which has the
exceptional divisor $E$ consisting of four $-1$ curves over the 4 points given by
$(0:\sqrt[4]{-1}:0:1)$. Consider the projection
$$p:S \setminus\{y_2 = x_1 = 0\} \hookrightarrow \bP^3 \setminus\{y_2 = x_1 = 0\}\to \bP^1 \ \ \mbox{given by} \ \ (y_2:y_1:x_1:x_0) \to (y_2:x_1).$$
Over $\{x_0 = 1\}$ one has an embedding of an open subset of $\tilde \bP^3$ into
$\bP^1 \times \A^3$, which yields an open embedding $e$ of an open subset $U$ of $\tilde S$ into
$\bP^1 \times \A^3$. Note that $\bP^1 \times \A^3$ is endowed with a natural projection $pr_1:
\bP^1 \times \A^3 \to \bP^1$. Over $U\setminus\{y_2= x_1 =0\}$ one has
$$p = pr_1 \circ e.$$
Thus by glueing, $p$ extends to a morphism $\tilde S \to \bP^1$, which is a family of projective curves of
degree 4. This family has a section $D = \{y_1=x_0=0\}$. One checks easily the singular loci of the
fibers do not meet $D$. By $\tilde S \times \F_3\to \bP^1$, we
have a family of surfaces. Let $\gamma_{\F_3}$ denote the generator of the Galois group of $\F_3 \to
\bP^1$, which acts via pullback by $\xi$ on $\omega_{\F_3}$. The quotient map onto
$\tilde S \times \F_3/(\gamma,\gamma_{\F_3})$ yields a quotient singularity of type $A_{3,2}$
(with the notation of \cite{Bart}). As in \cite{JCR} Section $9.2$ described one must
blow up $D$ and in a second step one blows up the fixed locus of the exceptional
divisor over $D$. By blowing down the image of the proper transform of the exceptional divisor over $D$,
one obtains a Calabi-Yau 3-manifold, which has obviously complex multiplication. 
\end{pkt}

\begin{pkt}
The automorphism $\gamma$ acts on
$\tilde S$, too. The quotient map $\varphi$ onto $M = \tilde S/\gamma$ is ramified over $E$ and
$D = \{y_1=x_0=0\}$. Since $D$ is a rational curve on a $K3$ surface, the adjunction formula implies
that $D.D=-2$. By the Hurwitz formula, one has
$$\varphi^*K_M \sim -2D-E.$$
Since
$$3\cdot K_M^2 = (\varphi^*K_M)^2,$$
one concludes that
$$c_1(M)^2 = K_M^2 = -4.$$
Thus the Noether formula tells us that $b_2(M) = 14$. Since we have blown up 4 points, one concludes
that $h^{1,1}_0(S) = 10$. Thus
$$h^{1,1}_1(S)=h^{1,1}_2(S) = 5.$$

By the fact that one has an exceptional divisor consisting of 12 copies of $\bP^2$ and 6 rational
ruled surfaces, one obtains similarly to \cite{JCR} Section $10.3$:
$$h^{1,1}(X) = h^{1,1}_0(S)+h^{1,1}_0(E) + 18 = 29$$
$$h^{2,1}(X) = h^{1,0}_1(E)\cdot h^{1,1}_2(S) = 5$$
\end{pkt}

\end{document}